\documentclass[11pt]{amsart}

\usepackage{amsfonts, amssymb, amscd}

\def\ZZ                 {{\mathbb Z}}
\def\QQ                 {{\mathbb Q}}

\def\CC                 {{\mathbb C}}

\def\mod			{{\,\rm mod}}
\def\ml				{{\mod l}}
\def\ii                 {{\rm i}}

\newcommand\draft[1]	{{}}

\newtheorem{lemma}{Lemma}[section]
\newtheorem{theorem}[lemma]{Theorem}
\newtheorem{corollary}[lemma]{Corollary} 
\newtheorem{proposition}[lemma]{Proposition}
\theoremstyle{definition}
\newtheorem{definition}[lemma]{Definition}

\newtheorem{remark}[lemma]{Remark}
\theoremstyle{remark}
\newtheorem*{proof*}{Proof}
\numberwithin{equation}{section}

\begin{document}
\title{On Wronskians of weight one Eisenstein series} 
 
\author{Lev A.  Borisov}
\address{Department of Mathematics, University of Wisconsin, 
417 Van Vleck Hall, 480 Lincoln Dr, Madison, WI 53706, USA.
{\tt email:borisov@math.wisc.edu}.}

\thanks{The author was partially supported by NSF grant DMS-0140172.}

\begin{abstract}
We describe the span of Hecke eigenforms of weight four with nonzero
central value of $L$-function in terms of Wronskians of certain weight
one Eisenstein series.
\end{abstract}

\maketitle
\section{Introduction}
For any positive integer $l$ we consider the congruence
subgroup $\Gamma_1(l)\subseteq Sl_2(\ZZ)$. The space of 
cusp forms for $\Gamma_1(l)$ 
of a given weight $k$ splits according to the eigenvalues
of Hecke operators.  We say that a Hecke eigenform has analytic
rank zero, if the central value of the corresponding 
$L$-function is nonzero.

It has been shown in \cite{vanish} that the span of Hecke eigenforms
of weight two coincides with the span of the cuspidal parts of products 
of certain weight $1$ Eisenstein series 
for the group $\Gamma_1(l)$. These series are the logarithmic derivatives in the $z$ direction of the standard $\theta$-function,
evaluated at $\frac al$ for $a=1,\ldots,l-1$. It is convenient to
look at the Fricke involutions of these Eisenstein series.
These are linear combinations of the original series and are
given by 
$$
s_a(q) = (\frac 12-\{\frac al\})+\sum_{n>0}q^n\sum_{d\vert n}
(\delta_d^{a\mod l}-\delta_d^{-a\mod l}),
$$
where $q=\exp(2\pi\ii\tau)$ and $\delta$ is a version of Kronecker
symbol.  In this paper we look at the Wronskians $W(s_a,s_b)$
defined as usual by $W(s_a,s_b)= (\frac d {d\tau}s_a)s_b
-(\frac d {d\tau}s_b)s_a$. It is easy to see that $W(s_a,s_b)$ 
is always a cusp form of weight $4$, and the main result of this
paper relates the span of such forms with the span of Hecke
eigenforms of analytic rank zero.

\smallskip
\noindent
{\bf Theorem \ref{main}.}
For arbitrary $l>1$ the span of Hecke eigenforms
of weight four and analytic rank zero is equal to the
span of the Wronskians $W(s_a(\tau),s_b(\tau))$ for all
$a,b\in \ZZ/l\ZZ$.

\smallskip

Before we explain the idea of the proof of this paper, 
we remark that it should be possible to prove Theorem \ref{main}
using Rankin-Selberg method, by combining
the formulas \cite[4.3, equation (4)]{Zagier} 
and \cite[Theorem 4.6.3]{Scholl}. However,
we chose to use the technique of \cite{vanish} and \cite{highweight}
that emphasizes the map from modular symbols to modular
forms. 

The space $M_4(l)$ 
of modular symbols of weight four can be thought of 
as a combinatorial counterpart to the space of modular forms.
It is a vector space of roughly twice the dimension, and it 
contains subspaces $S_4(l)_+$ and $S_4(l)_-$ which are 
naturally dual to the space ${\mathcal S}_4(l)$ of cusp forms
of weight four. Moreover, the action of Hecke operators 
on the space of modular symbols is given explicitly, see \cite{Merel}.
Ignoring minor complications due to old forms, 
the span of Hecke eigenforms of weight four and analytic rank zero
can be seen as the image of the endomorphism
$
\rho:{\mathcal S}_4(l)\to{\mathcal S}_4(l)
$
given by  
$$\rho(f)=\sum_{n>0} L(T_nf,2)q^n$$ where $T_n$ 
denote the Hecke operator. We observe that $L(f,2)q^n$ is 
the result of the pairing $\langle f,xy(0,1)_-\rangle$ of $f$
a certain element of $S_4(l)_-$ to calculate $\rho$ in terms 
of modular symbols as a composition of maps
$$
{\mathcal S}_4(l)
\stackrel {Int} \to (S_4(l)_-)^*
\stackrel {PD}\to 
S_4(l)_+\stackrel{\mu}\to {\mathcal S}_4(l)
$$
where $Int$ is induced by the integration
pairing of ${\mathcal S}_4(l)$ and $S_4(l)_-$,
the $PD$ is the Poincar\'e duality map which we 
define in Section \ref{sec.PD},
and $\mu$ is the Wronskian map, defined in Section \ref{sec.mu}.

The map $PD$ is a weight four analog of the intersection
pairing on weight $2$ symbols considered in \cite{vanish}. 
It is shown to be nondegenerate in Section \ref{sec.PD}
as a consequence of a modular symbol formula for Petersson
inner product. The map $\mu$ is the main novelty
of this paper. It is a map from the space 
of modular symbols to the space of modular forms, which in
particular maps $xy(a,b)$ to the Wronskian $W(s_a,s_b)$.
Our calculations are purely elementary and rely on properties
of the Euclid algorithm and some explicit calculations with
modular symbols. 

There are several directions in which one can try to 
extend the results of this paper. For example, one can look
at the subspaces in the spaces of modular forms of higher
weight that are spanned by Wronskians of Eisenstein series 
of higher weight. Intuition derived from \cite{highweight} 
and \cite{Zagier}
suggests that these would be related to values of the $L$-function
at $2$. Consequently, we expect the Wronskians to span the 
whole space in the higher weight setting. 

It is worth mentioning that the product and the Wronskian
are the first two cases of Cohen operators (see \cite{Zagier}).
One can wonder if the forms of higher weight of analytic
rank zero can be described in terms of higher Cohen operators
of $s_a$. Clearly, for a high enough weight this seems impossible
for dimension reasons. On the other hand, one could perhaps apply
Cohen operators to the theta function itself, rather than its logarithmic
derivatives, similar to the definition of $\mu$ on the noncuspidal
symbols of weight four. But this is all but a speculation at this
point. 

One might hope to use the construction of this paper to 
give upper bounds on the number of Hecke eigenforms of higher
analytic rank. However, analogous statements for weight two, 
at least so far, has not lead to such results. It can also be argued
that there may be some deeper reason behind the results of 
this paper and \cite{vanish} which is yet to be uncovered. From
this point of view, it would be tempting to try to see the 
sums  along the runs of Euclid algorithm as a calculation of 
an Euler characteristics of some complex, whose cohomology 
is located at top and bottom location only. But at the moment
we do not have a suitable candidate for it. Finally, one can wonder
whether derivatives of $L$-function at the central value can
be somehow seen in terms of Eisenstein series and Cohen
operators.

{\bf Notations.} We denote by $\mathcal H$ the upper half-plane
and denote by $\tau$, $\Im(\tau)>0$ the complex coordinate on it.
We use the notation
$q=\exp(2\pi\ii\tau)$ when writing Fourier expansion of 
modular forms.
Throughout the paper $l$ denotes the level,
and it is generally fixed, except for the proof of Theorem 
\ref{main} that requires induction on the level. We use
a slightly modified Kronecker $\delta$ notation 
$\delta_{u}^{v\mod w}$ which gives $1$ when 
$u=v\mod w$ and $0$ otherwise.

{\bf Acknowledgments.} This paper grew out of 
a search of a (weight two) 
skew-symmetric analog of \cite{vanish} which the author 
talked about on and off for a few years with Paul Gunnells.
The author also thanks Lo\"ic Merel for helpful remarks
regarding the Poincar\'e duality map.

\section{Modular symbols of weight four}
Our main reference for modular symbols is 
the paper \cite{Merel} by Merel,
which in turn builds on the work of Manin and Shokurov.
In this section we recall the purely combinatorial
description of modular (Manin, in the terminology 
of \cite{Merel}) symbols of 
weight four for the group $\Gamma_1(l)$. 

The modular symbols of weight four and level $l$ is 
a quotient of the vector space with basis $x^2(u,v)$,
$xy(u,v)$, $y^2(u,v)$, with $(u,v)\in (\ZZ/l\ZZ)^2$,
$gcd(u,v,l)=1$ by the span of the relations
\begin{equation}\label{rels}
\begin{array}{l}
x^2(u,v)+y^2(v,-u),~xy(u,v)-xy(v,-u),~
y^2(u,v)+x^2(v,-u)\\
xy(v,-u-v)-xy(-u-v,u)+y^2(-u-v,u)+x^2(u,v)-xy(u,v)
\end{array}
\end{equation}
for all $u,v$ with $gcd(u,v,l)=1$.

\begin{remark}
Our set of relations looks somewhat smaller than that of
\cite{Merel}, where the relations 
are 
$$
\begin{array}{l}P(x,y)(u,v)+P(y,-x)(u,v)\\
P(x,y)(u,v)+P(y-x,-x)(-u-v,u)+P(-y,x-y)(v,-u-v)
\end{array}
$$
for an arbitrary degree two homogeneous polynomial
$P(x,y)$. The "missing" relations are obtained by
cyclic permutations of $(u,v,-u-v)$ in the last line
of \eqref{rels},
so the two definitions of modular symbols are equivalent.
\end{remark}

Recall that the subspace $S_4(l)\subset M_4(l)$ 
of cuspidal modular symbols is characterized as
follows. The cusps of the modular curve $X_1(l)=
\overline{{\mathcal H}/\Gamma_1(l)}$
are in one-to-one correspondence with elements
of the set $I=\{(a,b),~a\in \ZZ/l\ZZ,~b\in (\ZZ/(a,l)\ZZ)^*\}/\pm$.
This correspondence maps $(a,b)$ to $\frac {b^*}a\in\QQ\cup{\ii\infty}$
where $b^{*}$ is the
inverse of $b\mod (a,l)$.
For every element of $I$ there is a map
$M_4(l)\to \CC$ defined by 
\begin{equation}\label{cuspmap}
\begin{array}{rcl}
x^2(u,v)&\mapsto& 
\delta_u^{a\mod l}\delta_v^{b\mod(a,l)}
+\delta_u^{-a\mod l}\delta_v^{-b\mod(a,l)}
\\
y^2(u,v)&\mapsto& -\delta_v^{a\mod l}\delta_u^{-b\mod(a,l)}
-\delta_v^{-a\mod l}\delta_u^{b\mod(a,l)}
\\
xy(u,v)&\mapsto& 0
\end{array}
\end{equation}
Then the space of cuspidal symbols $S_4(l)$ is
defined as the intersection of the kernels of all these maps. 

We are now ready to formulate the main result of this
section.
\begin{proposition}\label{xyspan}
The space of cuspidal symbols $S_4(l)$ is spanned by
the modular symbols of the form $xy(u,v)$.
\end{proposition}

\begin{proof}
We can use the first set of equations to solve for
$y^2(u,v)$. Then we can think of modular symbols of weight four
as being spanned by $x^2(u,v)$ and $xy(u,v)$,
subject to conditions
\begin{equation}\label{star}
\begin{array}{l}
x^2(u,v)-x^2(-u,-v)=0,~xy(u,v)-xy(v,-u)=0,~\\
x^2(u,u+v)-x^2(u,v)=xy(v,-u-v)-xy(-u-v,u)-xy(u,v).
\end{array}
\end{equation}

Clearly, $xy(u,v)$ are cuspidal. 
On the other hand,
if a linear combination of $w=\sum_{u,v}\alpha_{u,v}x^2(u,v)$ is cuspidal, then for each $a\mod l$ and each $b\mod (a,l)$
$$
\sum_{v=b\mod(a,l)} (\alpha_{u,v} +\alpha_{-u,-v}) = 0
$$
By using relations $x^2(u,v)=x^2(-u,-v)$ we can 
write $w$ as a linear combination of
$x^2(u,v+ku)-x^2(u,v+(k-1)u)$ which is then written 
as a linear combination of $xy(u',v')$.
\end{proof}

\begin{remark}\label{extrarels}
In addition to the obvious symmetry relations
$xy(u,v)=xy(v,-u)$ there are still some other linear relations 
among the symbols
$xy(u,v)$ in $S_4(l)$. In fact, one can show that for $l\geq 5$
the linear relations on $xy(u,v)$ in $M_4(l)$ (or $S_4(l)$)
are spanned by the symmetry relations and   
$$
\sum_{k=0}^{l-1} \Big(xy(v+ku,-(k+1)u-v)-xy(-(k+1)u-v,u)-xy(u,v+ku)\Big)=0
$$
for all $u$ and $v$ with  $gcd(u,v,l)=1$. We leave the proof of 
this claim to the reader, as it will not be used elsewhere in the 
paper.
\end{remark}

\draft{
\begin{proof}
First of all, it is clear that the above are indeed relations
on $xy(u,v)$. To show that every relation on $xy(u,v)$
is a linear combination of the above, note that it has to 
be a linear combination of the equations \eqref{star}
which has no $x^2(u,v)$ present. So one needs to study
linear combinations of $x^2(u,v)-x^2(-u,-v)$ and 
$x^2(u,u+v)-x^2(u,v)$ that are zero. Note that the problem
splits according to the elements of the set of cusps $I$.

Consider first the case of $\pm(u\mod l, v\mod (u,l))$
with $2u\neq 0\mod l$. Then we need to look at linear 
combinations of $x^2(u,w)-x^2(-u,-w)$,
$x^2(u,u+w)-x^2(u,w)$ and $x^2(-u,-u-w)-x^2(-u,-w)$
for $w=v\mod(u,l)$. We observe that the linear
combinations 
$$
(x^2(u,w)-x^2(-u,-w)) + (x^2(u,u+w)-x^2(u,w)) 
$$$$
+(x^2(-u,-u-w)-x^2(u,u+w))-(x^2(-u,-u-w)-x^2(-u,-w))
$$
lead to the following relations on $xy(u',v')$
$$
(xy(w,-u-w)-xy(-u-w,u)-xy(u,w))$$$$
-(xy(-w,u+w)-xy(u+w,-u)-xy(-u,-w))=0
$$
which are linear combinations of $xy(u_1,v_1)-xy(-v_1,u_1)=0$.
We can use these linear combinations to get rid of 
all $x^2(-u,-u-w)-x^2(-u,-w)$ terms, so we now only need
to find all linear relations between $x^2(u,w)-x^2(-u,-w)$
and $x^2(u,u+w)-x^2(u,w)$ for $w=v\mod (u,l)$. 
These now could only involve $x^2(u,u+w)-x^2(u,w)$,
and there is exactly one such relation up to scaling,
namely 
$$
\sum_{w=v\mod(u,l)} (x^2(u,u+w)-x^2(u,w)).
$$
We observe that a nonzero multiple of this relation yields
precisely
$$\sum_{k=0}^{l-1} \Big(xy(v+ku,-(k+1)u-v)-xy(-(k+1)u-v,u)-xy(u,v+ku)\Big)=0.$$

If $2u=0\mod l$, but $2v\neq 0\mod(u,l)$, then the argument
is essentially unchanged and is left to the reader. If $2v=0\mod(u,l)$ then, since $v$ is invertible $\mod (u,l)$,
we have $(u,l)=1$ or $(u,l)=2$. In either case $2(u,l)=0\mod l$,
so $l\vert 4$, which contradicts our assumption $l\geq 5$.
\end{proof}
}

We now recall that $M_4(l)$ and $S_4(l)$ naturally
split according to the eigenvalue of the involution
$i$ given by 
$$
x^2(u,v)\mapsto x^2(-u,v),~xy(u,v)\mapsto -xy(-u,v),~
y^2(u,v)\mapsto y^2(-u,v).
$$
We define the corresponding eigenspaces by $M_4(l)_+$,
$M_4(l)_-$,  $S_4(l)_+$ and $S_4(l)_-$. There are 
symmetrization maps $M_4(l)\to M_4(l)_\pm$ given
by $t\to \frac 12(t\pm i(t))$, and similarly for $S_4(l)\to S_4(l)_\pm$. We use a subscript to indicate the symmetrization
of a symbol.
We can now apply Proposition \ref{xyspan} to $S_4(l)_\pm$.
\begin{corollary}\label{genpm}
The space $S_4(l)_\pm$ is a linear span of the symbols
$xy(u,v)_\pm$ with $(u,v)\in (\ZZ/l\ZZ)^2$ and $gcd(u,v,l)=1$.
\end{corollary}

\begin{remark}
It is amusing to observe that for prime $l\geq 3$ the space
$S_4(l)_+$ is a linear span of symbols $xy(u,v)_+$ with 
$(u,v)\in (\ZZ/l\ZZ)^2-(0,0)$ with linear relations
among these symbols generated by 
$$
xy(u,v)_+=-xy(-u,v)_+=-xy(v,u)_+.
$$
Clearly, these relations hold in $S_4(l)_+$ and, by themselves,
they cut its dimension down to at most $\frac 18(l-1)(l-3)$.
On the other hand, by \cite{Merel}, $S_4(l)_+$
is dual to the space ${\mathcal S}_4(l)$
of cusp forms of weight four, which
has dimension $\frac 18(l-1)(l-3)$ by the usual 
Riemann-Roch calculation. This shows that all 
other relations on $xy(u,v)_+$ follow from the above 
symmetry relations (which can also be checked directly along
the lines of Remark \ref{extrarels}). 
One can thus identify $S_4(l)_+$
with the second exterior power
of the vector space of dimension $(l-1)/2$ which is 
generated by the symbols $r_a$ for $a\mod l$ with $r_{-a}=-r_a$.
\end{remark}

\section{Poincar\'{e} Duality for Modular Symbols}\label{sec.PD}
The goal of this section is to explicitly describe a certain 
map $PD:M_4(l)^*\to M_4(l)$ which is a weight four analog
of the Poincar\'e duality for the weight two cuspidal symbols.
It is rather easy to show that $PD(M_4(l))\subseteq S_4(l)$.
In fact, we will see that  $PD(M_4(l))=S_4(l)$, which is crucial
for the argument of this paper. This 
is proved by comparison of $PD$ and the expression of
the Petersson inner product of cusp forms of weight four in terms 
of their period integrals.

\begin{definition} The linear map $PD:M_4(l)^*\to M_4(l)$ 
is defined  by sending  any linear function $\phi:M_4(l)\to \CC$
to the element of $M_4(l)$ given by 
$$
\begin{array}{c}
\frac 1{24} \sum_{u,v\in \ZZ/l\ZZ}
(
\phi((y-x)^2(-v,u+v)-(y+x)^2(v,v-u))x^2(u,v)
\\
-2\phi( y(y-x)(-v,u+v)-(-y)(y+x)(v,v-u))xy(u,v)
\\ 
+\phi( y^2(-v,u+v)-(-y)^2(v,v-u)) y^2(u,v))
\end{array}
$$
where we adopt
the convention $\phi(P(x,y)(u,v))=0$ for any $P(x,y)$ if
$gcd(u,v,l)>1$.
\end{definition}

\begin{proposition}\label{skewsym}
The bilinear form on $M_4(l)$ induced by $PD$ is skew-symmetric.
Namely, for any $\phi,\lambda \in M_4(l)^*$ one has 
$$
\lambda(PD(\phi))=-\phi(PD(\lambda)).
$$
\end{proposition}

\begin{proof}
We can express $\lambda(PD(\phi))$ as
$$
\begin{array}{c}
\frac 1{24} \sum_{u,v\in \ZZ/l\ZZ}
(
\phi((y-x)^2(-v,u+v)-(y+x)^2(v,v-u))\lambda(x^2(u,v))
\\
-2\phi( y(y-x)(-v,u+v)-(-y)(y+x)(v,v-u))\lambda(xy(u,v))
\\ 
+\phi( y^2(-v,u+v)-(-y)^2(v,v-u)) \lambda(y^2(u,v))).
\end{array}
$$
We then use the relations $P(x,y)(c,d) = -P(y,-x)(-d,c)$ 
to rewrite the terms with $(v,v-u)$ in terms of $(u-v,v)$.
Afterwards, we use the relation 
$P(x,y)(c,d)=-P(y-x,-x)(-c-d,c)-P(-y,x-y)(d,-c-d)$ to further
rewrite them in terms of $(v,-u)$ and $(-u,u-v)$. Finally,
we use the relations $P(x,y)(c,d)=-P(y,-x)(d,-c)$ 
to write the result in terms of $(u,v)$ and $(u,u-v)$. 
We also rewrite the terms with $(-v,u+v)$ in terms of 
$(u,v)$ and $(-u-v,u)$. After simplifications, this gives
$$
\begin{array}{rl}
(y-x)^2(-v,u+v)-(y+x)^2(v,v-u) 
=&x^2(u-v,u)\\&-x^2(-u-v,u)
\\
y(y-x)(-v,u+v)-(-y)(y+x)(v,v-u) 
=&-x(y+x)(u-v,u)\\&-x(x-y)(-u-v,u)
\\
y^2(-v,u+v)-(-y)^2(v,v-u)
=&(y+x)^2(u-v,u)\\&-(x-y)^2(-u-v,u) 
\end{array}
$$
which allows us to write $\lambda(PD(\phi))$ as
$$
\begin{array}{c}
\frac 1{24}
\sum_{u,v}
(
\phi(
x^2(u-v,u)-x^2(-u-v,u)
)
\lambda(x^2(u,v))
\\
-2\phi( 
-x(y+x)(u-v,u)-x(x-y)(-u-v,u)
)\lambda(xy(u,v))
\\ 
+\phi( 
(y+x)^2(u-v,u)-(x-y)^2(-u-v,u) 
) \lambda(y^2(u,v))).
\end{array}
$$
It remains to switch the indexing in $\sum_{u,v}$ 
so that $\phi(...)$ becomes
$\phi(...(u,v))$ and simplify to get $-\phi(PD(\lambda))$.
\draft
{
$$
\begin{array}{c}
\frac 1{24}
\sum_{u,v}
(
\phi(x^2(u,v))
(\lambda(x^2(v,v-u))-\lambda(x^2(-v,u+v)))
\\
-2\phi(-x(y+x)(u,v))\lambda(xy(v,v-u))
+2\phi(x(x-y)(u,v))\lambda(xy(-v,u+v))
\\
+\phi((y+x)^2(u,v))\lambda(y^2(v,v-u)))
-\phi((x-y)^2(u,v))  \lambda(y^2(-v,u+v)))
\end{array}
$$
$$
\begin{array}{c}
=\frac 1{24}
\sum_{u,v}
(
\phi(x^2(u,v))
(\lambda((x+y)^2(v,v-u))-\lambda((x-y)^2(-v,u+v)))
\\
-2\phi(xy(u,v))(\lambda((-y)(x+y)(v,v-u))-\lambda((y-x)y(-v,u+v)))
\\
+\phi(y^2(u,v))(\lambda(y^2(v,v-u)))-\lambda(y^2(-v,u+v))))
\end{array}
$$
}
\end{proof}

\begin{proposition}\label{passthrough}
The map PD passes through the spaces of cusp symbols,
namely there is a commutative diagram
$$
\begin{array}{ccc}
M_4(l)^*&\stackrel{PD}\to&M_4(l)\\
\downarrow&&\uparrow\\
S_4(l)^*&\to&S_4(l)
\end{array}
$$
with the side maps coming from the natural inclusions 
$S_4(l)\to M_4(l)$.
\end{proposition}

\begin{proof}
First, let $\psi_{(a,b)}:M_4(l)\to \CC$ be the evaluation at the cusp
$\pm(a,b)$ with $a\in \ZZ/l\ZZ$ and $b\in \ZZ/(a,l)\ZZ$
as in \eqref{cuspmap}.
Let us show that 
$\psi_{(a,b)}(PD(\phi))=0$ for any $\phi$. Using \eqref{cuspmap},
we get 
$$
24\psi_{(a,b)}(PD(\phi))
=2\phi\Big(\sum_{v=b\mod(a,l)} (
(y-x)^2(-v,a+v)-(y+x)^2(v,v-a))
$$$$
-
\sum_{u=-b\mod(a,l)} (
y^2(-a,u+a)-y^2(a,a-u)
)
\Big)
$$
By switching from $u$ to $u\mp a$ in the two terms of the 
last sum, we see that it vanishes. For the first sum, we
switch from $v$ to $v-a$ in the first term. Then we use
$(x-y)^2(c,d)=-(x+y)^2(-d,c)$ to reduce it (up to a constant)
to the value of
$\phi$ on $\sum_{v=b\mod(a,l)}(y-x)^2(a-v,v)$. We now apply 
the relations on modular symbols to rewrite this as 
$$
-\sum_{v=b\mod(a,l)} (x^2(v,-a)  + y^2(-a,a-v))
$$$$
=-\sum_{v=b\mod(a,l)} (x^2(v,-a)  + y^2(a,v))=0.
$$
This shows that the image of $PD$ sits inside $S_4(l)$.

By Proposition \ref{skewsym}, we now see that 
$\phi(PD(\psi_{a,b}))=0$ for any $\phi$. 
This shows that $PD$ passes through
$S_4(l)^*$ which finishes the proof.
\end{proof}

The key result of this section hinges on a formula
for the  Petersson inner product
of cusp forms in terms of period integrals.
Recall that the Petersson inner product of two holomorphic cusp forms of weight four with respect to 
$\Gamma_1(l)$ is defined as 
$$(f,g)_{\rm Petersson} = 
\int\hskip-.2cm\int_{{\mathcal H}/\Gamma_1(l)}
f(\tau)\overline{g(\tau)} \Im(\tau)^2d\tau d\bar\tau.
$$
Period integrals define
a pairing between the space ${\mathcal S}_4(l)$
of cusp forms of weight four
and $M_4(l)$, which we denote by $\langle~,~\rangle$. 
This pairing is a crucial feature of the theory of modular 
symbols, and we refer the reader to \cite{Merel} for its
definitions and properties.
\begin{theorem}\label{Petmod}
For any two holomorphic weight four forms $f$ and $g$ 
there holds
$$
\begin{array}{rl}
{\rm (1)~}&(f,g)_{\rm Petersson}
=-\frac 1{24}\sum_{c,d\in \ZZ/l\ZZ,~gcd(c,d,l)=1} 
\\ &
\Big(
\overline{(
\langle g, (y-x)^2(-d,c+d)\rangle - 
\langle g, (y+x)^2(d,d-c)\rangle
)}
\langle f,x^2(c,d)\rangle
\\ &
-2\overline{(
\langle g, y(y-x)(-d,c+d)\rangle - 
\langle g, -y(y+x)(d,d-c)\rangle
)}
\langle f,xy(c,d)\rangle
\\ &
+\overline{(
\langle g, y^2(-d,c+d)\rangle - 
\langle g, (-y)^2(d,d-c)\rangle
)}
\langle f,y^2(c,d)\rangle
\Big)
\end{array}
$$
$$
\begin{array}{rl}
{\rm (2)~}&0
=-\frac 1{24}\sum_{c,d\in \ZZ/l\ZZ,~gcd(c,d,l)=1} 
\\ &
\Big(
{(
\langle g, (y-x)^2(-d,c+d)\rangle - 
\langle g, (y+x)^2(d,d-c)\rangle
)}
\langle f,x^2(c,d)\rangle
\\ &
-2{(
\langle g, y(y-x)(-d,c+d)\rangle - 
\langle g, -y(y+x)(d,d-c)\rangle
)}
\langle f,xy(c,d)\rangle
\\ &
+{(
\langle g, y^2(-d,c+d)\rangle - 
\langle g, (-y)^2(d,d-c)\rangle
)}
\langle f,y^2(c,d)\rangle
\Big).
\end{array}
$$
\end{theorem}

\begin{proof}
Consider the cosets $\Gamma_1(l)\lambda$ for $\lambda\in
Sl_2(\ZZ)$. Then the fundamental domain 
${\mathcal H}/\Gamma_1(l)$ can be chosen as 
$\bigcup_{\Gamma_1(l)\lambda}\lambda(D_0)$ for any fundamental
domain $D_0$ of $\Gamma_1(l)$. Moreover, we can use a 
union of three different such $D_0$ to write the Petersson 
pairing as
$$
(f,g)_{\rm Petersson}=\frac 13\sum_{\Gamma_1(l)\lambda}
\int\hskip-.2cm\int_{\lambda(D)}
f(\tau)\overline{g(\tau)} \Im(\tau)^2d\tau d\bar\tau
$$
where $D$ is the geodesic triangle in $\mathcal H$ with 
vertices $\ii\infty$, $-1$ and $0$. The boundary
of $D$ consists of the vertical lines $\Re(\tau)=0$ and
$\Re(\tau)=-1$, as well as the upper half of the circle of 
radius $\frac 12$ centered at $-\frac 12$.

We use $\Im(\lambda(\tau))^2=\Im(\tau)^2\vert c\tau+d\vert^{-4}$
where $\lambda(\tau) =\frac {a\tau+b}{c\tau+ d}$ to rewrite 
each term of the above sum as 
$$
\int\hskip-.2cm\int_{D}
f(\lambda(\tau))\overline{g(\lambda(\tau))} \Im(\tau)^2\vert c\tau+d\vert^{-8} d\tau d\bar\tau.
$$
For each such $\lambda$ we introduce for $i=0,1,2$
$$
G_{i,\lambda}(\tau)=\overline{\int_{-1}^{\tau} 
g(\lambda(s))(cs+d)^{-4}s^i ds}
$$
and $f_{i,\lambda}(\tau))=f(\lambda(\tau))(c\tau+d)^{-4}\tau^i$.
Then we write $\Im(\tau)^2 = -\frac 14 (\tau-\bar\tau)^2$ and 
use Stokes's Theorem to derive
$$
(f,g)_{\rm Petersson}=-\frac 1{12}\sum_{\Gamma_1(l)\lambda}
\int_{\partial D} \Big(G_{0,\lambda}(\tau)f_{2,\lambda}(\tau)
-2G_{1,\lambda}(\tau)f_{1,\lambda}(\tau)
$$$$+G_{2,\lambda}(\tau)f_{0,\lambda}(\tau)\Big)d\tau.
$$

The boundary of $D$ splits into three geodesics, and
our first claim is that the terms of the integration for
the $\int_{\ii\infty}^{-1}$ and $\int_{-1}^0$ of $\partial D$ 
cancel each other.
Consider the map $\sigma(\tau):=-\frac {1}{\tau+2}$.
Element $\sigma\in Sl_2(\ZZ)$ 
acts on the set of cosets $\Gamma_1(l)\lambda$ by right multiplication. We
observe that 
$$
\int_{-1}^0(G_{0,\lambda\sigma}(\tau)f_{2,\lambda\sigma}(\tau)
-2G_{1,\lambda\sigma}(\tau)f_{1,\lambda\sigma}(\tau)
+G_{2,\lambda\sigma}(\tau)f_{0,\lambda\sigma}(\tau))
d\tau
$$
$$
=\int_{-1}^0(G_{0,\lambda}(\sigma(\tau))f_{2,\lambda}(\sigma(\tau))
-2G_{1,\lambda}(\sigma(\tau))f_{1,\lambda}(\sigma(\tau))
$$$$+G_{2,\lambda}(\sigma(\tau))f_{0,\lambda}(\sigma(\tau)))
d\sigma(\tau)
$$
$$
=\int_{-1}^{\ii\infty}(G_{0,\lambda}(\tau)f_{2,\lambda}(\tau)
-2G_{1,\lambda}(\tau)f_{1,\lambda}(\tau)
+G_{2,\lambda}(\tau)f_{0,\lambda}(\tau))
d\tau.
$$
\draft
{
$$
G_{i,\lambda\sigma}(\tau)=\overline
{\int_{-1}^{\tau}
g(\lambda(-\frac 1{s+2}))(ds-c+2d)^{-4}s^ids
}
$$
$$
=
\overline
{\int_{-1}^{\sigma(\tau)}
g(\lambda(s))(cs+d)^{-4}s^{2-i}(-1-2s)^ids
}
$$
$$
G_{0,\lambda\sigma}(\tau) = 
G_{2,\lambda}(\sigma(\tau))
,~
G_{1,\lambda\sigma}(\tau) = 
-G_{1,\lambda}(\sigma(\tau))
-2G_{2,\lambda}(\sigma(\tau))
$$
$$
G_{2,\lambda\sigma}(\tau) = 
G_{0,\lambda}(\sigma(\tau))
+4G_{1,\lambda}(\sigma(\tau))
+4G_{2,\lambda}(\sigma(\tau))
$$
$$
f_{0,\lambda\sigma}(\tau)d\tau = 
f_{2,\lambda}(\sigma(\tau))d\sigma(\tau)
$$$$
f_{1,\lambda\sigma}(\tau)d\tau = 
(-f_{1,\lambda}(\sigma(\tau))
-2f_{2,\lambda}(\sigma(\tau)))d\sigma(\tau)
$$
$$
f_{2,\lambda\sigma}(\tau)d\tau  = 
(f_{0,\lambda}(\sigma(\tau))
+4f_{1,\lambda}(\sigma(\tau))
+4f_{2,\lambda}(\sigma(\tau)))d\sigma(\tau)
$$
}

\noindent
The first equality is verified by a lengthy but straightforward 
calculation, which is left to the reader, since we will perform
a similar calculation below. It is crucial that
we chose $(-1)$ as
the lower limit of integration in the definition of $G_{i,\lambda}$
and that $\sigma$ preserves $(-1)$.

So now we are left with
$$
(f,g)_{\rm Petersson}=-\frac 1{12}\sum_{\Gamma_1(l)\lambda}
\int_{0}^{\ii\infty} \Big(G_{0,\lambda}(\tau)f_{2,\lambda}(\tau)
-2G_{1,\lambda}(\tau)f_{1,\lambda}(\tau)
$$
$$
+G_{2,\lambda}(\tau)f_{0,\lambda}(\tau)\Big)d\tau.
$$
We will do a similar trick, this time with $\nu(\tau)=-\frac 1\tau$
instead of $\sigma(\tau)$. It has an effect of switching the 
direction of integration, but since it does not preserve $(-1)$,
the functions $G_{i,\lambda}$ acquire extra additive terms.
More specifically, one has
$$
G_{i,\lambda\nu}(\tau)=\hskip -2pt
\overline
{\int_{-1}^{\tau}
g(\lambda(-\frac 1{s}))(ds-c)^{-4}s^ids
}
=\hskip -2pt
\overline
{\int_{1}^{\nu(\tau)}\hskip -.2cm
g(\lambda(s))(cs+d)^{-4}s^{2-i}(-1)^ids
}
$$
$$
=
(-1)^iG_{2-i,\lambda}(\nu(\tau))
+
\overline
{\int_{1}^{-1}
g(\lambda(s))(cs+d)^{-4}s^{2-i}(-1)^ids
}
$$
$$
f_{i,\lambda\nu}(\tau)d\tau = (-1)^i
f_{2-i,\lambda}(\nu(\tau))d\nu(\tau).
$$
We rewrite $(f,g)_{\rm Petersson}$ as
$$
-\frac 1{24}\Big (\sum_{\Gamma_1(l)\lambda}
\int_{0}^{\ii\infty} \Big(G_{0,\lambda}(\tau)f_{2,\lambda}(\tau)
-2G_{1,\lambda}(\tau)f_{1,\lambda}(\tau)
+G_{2,\lambda}(\tau)f_{0,\lambda}(\tau)\Big)d\tau
$$$$+\sum_{\Gamma_1(l)\lambda}
\int_{0}^{\ii\infty} \Big(G_{0,\lambda\nu}(\tau)f_{2,\lambda\nu}(\tau)
-2G_{1,\lambda\nu}(\tau)f_{1,\lambda\nu}(\tau)
+G_{2,\lambda\nu}(\tau)f_{0,\lambda\nu}(\tau)\Big)d\tau\Big)
$$
which together with transformation formulas for $G$ and $f$
implies, after cancelling $\pm\int_{0}^{\ii\infty}$
$$
(f,g)_{\rm Petersson}=-\frac 1{24}
\sum_{\Gamma_1(l)\lambda}
\Big(
\overline{
\int_{1}^{-1} g(\lambda(s)) (cs+d)^{-4} s^2 ds
}
\int_{0}^{\ii\infty}f_{0,\lambda}(\tau)d\tau
$$
$$
-2\overline{
\int_{1}^{-1} g(\lambda(s)) (cs+d)^{-4} s ds
}
\int_{0}^{\ii\infty}f_{1,\lambda}(\tau)d\tau
$$
$$
+\overline{
\int_{1}^{-1} g(\lambda(s)) (cs+d)^{-4} ds
}
\int_{0}^{\ii\infty}f_{2,\lambda}(\tau)d\tau
\Big).
$$
It remains to write $\int_{1}^{-1}$ in terms of pairings with 
modular symbols by writing the arc from $1$ to $(-1)$ in terms 
of the unimodular arcs from $1$ to $\ii\infty$ and from $\ii\infty$
to $(-1)$. Namely, for a homogeneous 
degree two polynomial $P(x,y)$ one has
$$
\int_{1}^{-1} g(\lambda(s))(cs+d)^{-4}P(s,1)ds
$$$$
=
\int_{1}^{0}g(\lambda(s))(cs+d)^{-4}P(s,1)ds
-\int_{-1}^{0}g(\lambda(s))(cs+d)^{-4}P(s,1)ds
$$
$$
=
\int_{0}^{\ii\infty}g(\lambda(\frac 1{1-s}))(-ds+ (c+d))^{-4}
P(1,1-s)ds
$$
$$
-\int_{0}^{\ii\infty}g(\lambda(-\frac 1{1+s}))
(ds+(d-c))^{-4}P(-1,1+s)dt
$$
$$
=\langle g, P(y,y-x)(-d,c+d)\rangle 
-\langle g, P(-y,y+x)(d,d-c) \rangle.
$$
Finally, one observes that cosets $\Gamma_1(l)\lambda$
are in one-to-one correspondence with pairs $(c,d)$ 
with $gcd(c,d,l)=1$, and the first claim of the theorem follows.

The second claim of the theorem is proved by the same technique.
This time we define $G_{i,\lambda}(\tau)$ as $\int_{-1}^\tau
g(\lambda(s))(cs+d)^{-4}s^ids$. Consequently, it is holomorphic,
and the Stokes's Theorem gives $0$ instead of the Petersson
product. The rest of the calculations are unchanged.
\end{proof}

\begin{remark}
Similar formulas for Petersson product 
are already present in the literature. They seem to
go back to at least as far as \cite{Haberland}
and \cite{KohnenZagier}. We learned the argument 
(in weight two case) from \cite{Merel2001}. In addition
to extending it to weight four, we 
streamlined it just slightly by looking 
at the union of three fundamental domains for $Sl_2(\ZZ)$, 
rather than one. This allowed us to avoid
integration between elliptic points.
\end{remark}

\begin{corollary}\label{sur}
The pairing on $S_4(l)$ induced by $PD$ is nondegenerate.
\end{corollary}

\begin{proof}
By \cite{Shokurov}, the integration pairing is a perfect pairing
between $S_4(l)$ and the direct sum $V_{hol}\oplus \bar V_{hol}$
of the spaces of holomorphic and anti-holomorphic cusp forms.
Every element  $\phi\in S_4(l)^*$ can therefore be written as
$\langle f, \cdot \rangle + \langle \bar g ,\cdot \rangle$. Suppose 
$PD(\phi)=0$. Denote by $\bar~$ the anti-isomorphism of
$M_4(l)$ that sends $\alpha x^iy^{2-i}(u,v)$ to $\bar\alpha x^i
y^{2-i}(u,v)$. Then Theorem \ref{Petmod} shows that 
$$
0=\phi(\overline{PD(\phi)})=
-\langle f,f\rangle_{\rm Petersson}
-\langle g,g\rangle_{\rm Petersson}\leq 0
$$
with equality only if $f=g=0$.
\end{proof}

\begin{remark}
The arguments of our proof 
of Theorem \ref{Petmod} extend naturally to arbitrary
integer weights $k\geq 2$
and arbitrary subgroups $\Gamma$ of finite index in $Sl_2(\ZZ)$.
We expect Corollary \ref{sur} to extend to arbitrary weight
and to arbitrary group $\Gamma$, after an appropriate definition
of $PD$. One would need to interpret the arguments 
of Propositions \ref{skewsym} and \ref{passthrough} to extend
them to this more general setting. For instance, we expect that 
the bilinear form on 
$S_k(\Gamma)$ induced by $PD$ is $(-1)^{k+1}$-symmetric.
We believe that maps $PD$ can be interpreted as an intersection
pairings in the middle cohomology of the Kuga varieties, although
we do not need this for the purposes of this paper. Nevertheless,
this is why we refer to $PD$ as the Poincar\'e duality map.
\end{remark}

The map $PD$ behaves well with respect to the involution 
$i$. We denote by $M_4(l)_\pm^*$ the eigenspaces of 
the dual involution $i^*$ on $M_4(l)^*$.
\begin{proposition}\label{pm}
$PD(M_4(l)_\pm^*)\subseteq M_4(l)_{\mp}$.
\end{proposition}

\begin{proof}
If $\phi\in M_4(l)_\pm^*$, then for any modular symbol
$P(x,y)(u.v)$ there holds $\phi(P(x,y)(u,v))=\phi(P(x,y)(u,v)_\pm)$.
Consequently,
$$
PD(\phi)=\frac 1{24}\sum_{u,v\in \ZZ/l\ZZ}(
\phi((y-x)^2(-v,u+v)_\pm-(y+x)^2(v,v-u)_\pm)x^2(u,v)
$$$$-2\phi( y(y-x)(-v,u+v)_\pm-(-y)(y+x)(v,v-u)_\pm)xy(u,v)
$$$$+\phi( y^2(-v,u+v)_\pm-(-y)^2(v,v-u)_\pm) y^2(u,v))
$$
$$
=\frac 1{24}\sum_{u,v\in \ZZ/l\ZZ}(
\phi((y-x)^2(-v,u+v)_\pm\mp(y-x)^2(-v,v-u)_\pm)x^2(u,v)
$$$$-2\phi( y(y-x)(-v,u+v)_\pm\mp(-y)(y-x)(-v,v-u)_\pm)xy(u,v)
$$$$+\phi( y^2(-v,u+v)_\pm\mp(-y)^2(-v,v-u)_\pm) y^2(u,v))
$$
$$
=\frac 1{24}
\sum_{u,v\in \ZZ/l\ZZ}(
\phi((y-x)^2(-v,u+v)_\pm)(x^2(u,v)\mp x^2(-u,v))
$$$$-2\phi( y(y-x)(-v,u+v)_\pm)(xy(u,v)\pm xy(-u,v))
$$$$+\phi( y^2(-v,u+v)_\pm) (y^2(u,v)\mp y^2(-u,v)))
$$
$$
=\frac 1{12}
\sum_{u,v\in \ZZ/l\ZZ}(
\phi((y-x)^2(-v,u+v)_\pm)x^2(u,v)_\mp
$$$$-2\phi( y(y-x)(-v,u+v)_\pm)xy(u,v)_\mp
+\phi( y^2(-v,u+v)_\pm) y^2(u,v)_\mp).
$$
\end{proof}

\begin{remark}
In what follows, we will abuse the notations somewhat 
to denote the induced map $M_4(l)_-^*\to M_4(l)_+$
by $PD$ as well. By Proposition \ref{passthrough},
this map comes from a map $S_4(l)_-^*\to S_4(l)_+$.
\end{remark}

\begin{corollary}\label{surj}
The induced map $PD:S_4(l)_-^*\to S_4(l)_+$ 
is an isomorphism.
\end{corollary}

\begin{proof}
Combine Corollary \ref{sur} and Proposition \ref{pm}.
\end{proof}

\section{The Wronskian map}\label{sec.mu}
In this section we define the map from the modular symbols
of weight four to cusp forms of weight four. 
First, we need to define the Eisenstein 
series $s_a(q)$, $t_a(q)$ and $r_a(q)$ for $a\in\ZZ/l\ZZ$.
Our notation for $s_a$ 
differs from that of \cite{vanish} by a Fricke 
involution. We recall that quasimodular forms of weight two
are linear combinations of the usual modular forms of weight two
and the (slightly non-modular) Eisenstein series $E_2(q)$
of weight $2$. In weight one, quasimodular forms are modular.

\begin{proposition}\label{rst}
For each $a\mod l$
there exist
$\Gamma_1(l)$-quasimodular forms $s_a(q)$, $t_a(q)$ and $r_a(q)$ 
of weights $1$, $2$ and $2$ respectively
given by 
$$
s_a(q) = (\frac 12-\{\frac al\})+\sum_{n>0}q^n\sum_{d\vert n}
(\delta_d^{a\mod l}-\delta_d^{-a\mod l}),
$$$$~{\rm if}~a\neq 0\ml,
s_0(q)=0,
$$
$$
t_a(q) =
constant
+\sum_nq^n\sum_{d\vert n}
\frac nk(\delta_{d}^{a\ml}+\delta_d^{-a\ml})
$$
$$
r_a(q) =  
constant+
\sum_nq^n\sum_{d\vert n}
d(\delta_{d}^{a\mod l}+\delta_d^{-a\mod l}),
$$
where the exact values of the constants depend on $a$
and $l$ and are determined uniquely by the quasimodularity.
\end{proposition}

\begin{proof}
These series are obtained as linear combinations of 
the weight one and two Eisenstein series considered
in \cite{vanish}. Details are left to the reader.
\end{proof}

\begin{definition}\label{mapmu}
We define the map 
$\mu:S_4(l)\to {\mathcal S}_4(l)$
by the formula
$$
x^2(u,v)\mapsto -2t_ur_v-\frac 1l q\frac {\partial}{\partial q}r_{v}
-\delta_v^{0\ml}q\frac {\partial}{\partial q}t_u,~
xy(u,v)\mapsto \frac 1{2\pi\ii}W(s_u,s_v),~
$$$$
y^2(u,v)\mapsto 2r_ut_v+\frac 1l q\frac {\partial}{\partial q}r_{u}
+\delta_u^{0\ml}q\frac {\partial}{\partial q}t_v.
$$
Note that the ring of quasimodular forms is closed
under $q\frac{\partial}{\partial q}$.
We will show in Theorem \ref{muisOK} below 
that $\mu$ is well-defined, i.e. it is compatible 
with the relations on modular symbols.
\end{definition}

\begin{theorem}\label{muisOK}
The map $\mu$ of definition \ref{mapmu} is well-defined.
\end{theorem}

\begin{proof}
We need to check that $\mu$ maps the relations \eqref{rels}
to zero. We have
$$
\mu(x^2(u,v)+y^2(v,-u))= 0 
$$
and 
$$
\mu(xy(u,v)-xy(v,-u))=\frac 1{2\pi\ii}(W(s_u,s_v)-W(s_v,s_{-u}))=0
$$
by the symmetry properties $s_{-a}=-s_a$, $r_a=r_{-a}$ 
and $t_a=t_{-a}$.

The difficult part is to show that $\mu$ maps
$$R=xy(v,-u-v)-xy(-u-v,u)+y^2(-u-v,u)+x^2(u,v)-xy(u,v)
$$
to zero. For each positive integer $n$ let us denote by $I(n)$ the
set of fourtuples $(m_1,k_1,m_2,k_2)\in \ZZ_{>0}^4$
that satisfy $m_1k_1+m_2k_2=n$. Let us denote
by $\sim$ the equality of power 
series in $q$ up to linear combinations 
of quasimodular forms of weights $0$, $1$, $2$, and the 
derivatives
of $s_a(q)$ with respect to $\tau$. This allows us to 
avoid looking at the specific values of the constant terms in
Proposition \ref{rst}. 

We have
$$\mu(R)\sim
-\frac 1l q\frac {\partial}{\partial q}r_{v}
-\delta_v^{0\ml}q\frac {\partial}{\partial q}t_u
+\delta_{u+v}^{0\ml}q\frac{\partial}{\partial q}t_u
+\frac 1l q\frac {\partial}{\partial q}r_{u+v}
$$$$
+\sum_{n>0}q^n
\sum_{I_n}
\Big(
(m_1k_1-m_2k_2)(\delta_{k_1}^{v\ml}-\delta_{k_1}^{-v\ml})
(\delta_{k_2}^{-u-v\ml}-\delta_{k_2}^{u+v\ml})
$$$$
-(m_1k_1-m_2k_2)(\delta_{k_1}^{-u-v\ml}-\delta_{k_1}^{u+v\ml})
(\delta_{k_2}^{u\ml}-\delta_{k_2}^{-u\ml})
$$$$
+2k_1m_2(\delta_{k_1}^{-u-v\ml}+\delta_{k_1}^{u+v\ml})
(\delta_{k_2}^{u\ml}+\delta_{k_2}^{-u\ml})
$$$$
-2m_1k_2(\delta_{k_1}^{u\ml}+\delta_{k_1}^{-u\ml})
(\delta_{k_2}^{v\ml}+\delta_{k_2}^{-v\ml})
$$
$$
-(m_1k_1-m_2k_2)(\delta_{k_1}^{u\ml}-\delta_{k_1}^{-u\ml})
(\delta_{k_2}^{v\ml}-\delta_{k_2}^{-v\ml})
\Big).
$$
We introduce the notation $A_{k_1,k_2}=
\delta_{k_1}^{u\ml}\delta{k_2}^{v\ml}+\delta_{k_1}^{-u\ml}
\delta_{k_2}^{-v\ml}$ to rewrite the above as
$$
\mu(R)\sim 
-\frac 1l q\frac {\partial}{\partial q}r_{v}
-\delta_v^{0\ml}q\frac {\partial}{\partial q}t_u
+\delta_{u+v}^{0\ml}q\frac{\partial}{\partial q}t_u
+\frac 1l q\frac {\partial}{\partial q}r_{u+v}
$$$$
+\sum_{n>0}q^n\sum_{I(n)}
\Big(
(m_2k_2-m_1k_1-2m_1k_2)A_{k_1,k_2}
$$$$
+(m_2k_2-m_1k_1+2k_1m_2)A_{k_2,-k_1-k_2}
+(m_1k_1-m_2k_2)A_{-k_1-k_2,k_1}
$$$$
-(m_2k_2-m_1k_1+2m_1k_2)A_{-k_1,k_2}
-(m_2k_2-m_1k_1-2k_1m_2)A_{k_2,k_1-k_2}
$$$$-(m_1k_1-m_2k_2)A_{k_1-k_2,-k_1}
\Big)
$$
We recall (see \cite{highweight}) that $I(n)$ is 
a disjoint union of the runs of Euclid algorithm.
The algorithm is given by the partially defined
map $up:I(n)\to I(n)$ 
$$
up:(m_1,k_1,m_2,k_2)\mapsto\Big\{
\begin{array}{ll}
(m_2,k_1+k_2,m_1-m_2,k_1),&m_1>m_2\\
(m_2-m_1,k_2,m_1,k_1+k_2),&m_1<m_2.
\end{array}
$$
Repeated applications of this map go from 
the subset of $I(n)$ with $k_1=k_2$ to the subset
of $I(n)$ with $m_1=m_2$, where $up$ is not defined.
As in \cite{highweight}, we will show that 
for each run of the algorithm the above sum is telescoping.
Namely, the "plus"
terms with $A_{k_1,k_2}$, $A_{k_2,-k_1-k_2}$,
$A_{-k_1-k_2,k_1}$ for $(m_1,k_1,m_2,k_2)$ cancel
the "minus"
terms with $A_{-k_1,k_2}$, $A_{k_2,k_1-k_2}$,
$A_{k_1-k_2,-k_1}$ for $up(m_1,k_1,m_2,k_2)$.
There are two cases to check, depending on whether 
$m_1>m_2$ or $m_1<m_2$. In the case of $m_1>m_2$
the "minus" terms for 
$up(m_1,k_1,m_2,k_2)=(m_2,k_1+k_2,m_1-m_2,k_1)$
are 
$$
-((m_1-m_2)k_1-m_2(k_1+k_2)+2m_2k_1)A_{-k_1-k_2,k_1}
$$$$
-((m_1-m_2)k_1-m_2(k_1+k_2)-2(k_1+k_2)(m_1-m_2))A_{k_1,k_2}
$$
$$
-(m_2k_1+m_2k_2)-(m_1-m_2)k_1)A_{k_2,-k_1-k_2}
$$
$$
=-(m_1k_1-m_2k_2)A_{-k_1-k_2,k_1}
-(-k_1m_1+k_2m_2-2k_2m_1)A_{k_1,k_2}
$$
$$
-(m_2k_2-m_1k_1+2m_2k_1)A_{k_2,-k_1-k_2}
$$
which is seen to equal the "plus" terms for $(m_1,k_1,m_2,k_2)$.
The case of $m_1<m_2$ is similar and left to the reader.
One needs to use there the symmetry 
$A_{k_1,k_2}=A_{-k_1,-k_2}$.
Consequently, the only terms that will not be cancelled are
the "plus" terms for the subset of $I(n)$ with $m_1=m_2$
and the "minus" terms for the subset of $I(n)$ with $k_1=k_2$.
This gives 
$$
\mu(R)\sim 
-\frac 1l q\frac {\partial}{\partial q}r_{v}
-\delta_v^{0\ml}q\frac {\partial}{\partial q}t_u
+\delta_{u+v}^{0\ml}q\frac{\partial}{\partial q}t_u
+\frac 1l q\frac {\partial}{\partial q}r_{u+v}
$$
$$
+\sum_{n>0}q^n
\Big(
\sum_{\substack{m_1,k_1,k_2>0\\m_1(k_1+k_2)=n}}(-nA_{k_1,k_2}
+nA_{k_2,-k_1-k_2}
+m_1(k_1-k_2)A_{-k_1-k_2,k_1})
$$$$
-
\sum_{\substack{m_1,m_2,k_1>0\\(m_1+m_2)k_1=n}}
(nA_{-k_1,k_1}
-nA_{k_1,0}
+(m_1-m_2)k_1A_{0,-k_1})
\Big)
$$
$$
=-\frac 1l q\frac {\partial}{\partial q}r_{v}
-\delta_v^{0\ml}q\frac {\partial}{\partial q}t_u
+\delta_{u+v}^{0\ml}q\frac{\partial}{\partial q}t_u
+\frac 1l q\frac {\partial}{\partial q}r_{u+v}
$$$$
+\sum_{n>0}nq^n
\sum_{d\vert n}
\Big(
\sum_{0<k\leq d}
(-A_{k,d-k}
+A_{d-k,-d})
-A_{0,-d}
$$$$+
\sum_{0<k\leq d}
(\frac{2k}d-1)A_{-d,k}
+\frac nd(A_{d,0}
-A_{-d,d})
\Big).
$$

We observe that 
$$
\sum_{0<k\leq d}\delta_k^{u\ml}=
\frac dl - \{\frac{d-u}l\} +\{-\frac ul\}
$$
and 
$$
\sum_{0<k\leq d}(\frac {2k}d-1)\delta_k^{u\ml}
=\frac ld\Big(\{\frac {d-u}l\}^2-\{\frac {d-u}l\}
-\{-\frac {u}l\}^2+\{-\frac ul\}\Big)
$$$$+\Big(1-\{\frac {d-u}l\}
-\{-\frac ul\}\Big).
$$
Consequently, 
$$
\mu(R)\sim
-\frac 1l q\frac {\partial}{\partial q}r_{v}
-\delta_v^{0\ml}q\frac {\partial}{\partial q}t_u
+\delta_{u+v}^{0\ml}q\frac{\partial}{\partial q}t_u
+\frac 1l q\frac {\partial}{\partial q}r_{u+v}
$$$$
+\sum_{n>0}nq^n
\sum_{d\vert n}
\Big(-\delta_{u}^{0\ml}(\delta_{d}^{v\ml}+\delta_d^{-v\ml})
$$$$
-\delta_{d}^{u+v\ml}
(\frac dl - \{\frac vl\} +\{-\frac ul\})
-\delta_{d}^{-u-v\ml}
(\frac dl - \{-\frac vl\} +\{\frac ul\})
$$$$
+\delta_{d}^{v\ml}
(\frac dl - \{-\frac{u}l\} +\{-\frac {u+v}l\})
+\delta_{d}^{-v\ml}
(\frac dl - \{\frac{u}l\} +\{\frac {u+v}l\})
$$$$+
\frac nd\delta_{v}^{0\ml}(\delta_{d}^{u\ml}+\delta_d^{-u\ml})
-\frac nd \delta_{u+v}^{0\ml}(\delta_d^{v\ml}+\delta_d^{-v\ml})
+ 
\delta_d^{-u\ml}
\cdot$$$$\cdot
\Big(\frac ld(\{\frac {-u-v}l\}^2-\{\frac {-u-v}l\}
-\{-\frac {v}l\}^2+\{-\frac vl\})
+(1-\{\frac {-u-v}l\}
-\{-\frac vl\})\Big)
$$
$$
+\delta_d^{u\ml}
\Big(\frac ld(\{\frac {u+v}l\}^2-\{\frac {u+v}l\}
-\{\frac {v}l\}^2+\{\frac vl\})
+(1-\{\frac {u+v}l\}
-\{\frac vl\})\Big)\Big).
$$
We use $\{t\}+\{-t\}=1-\delta_{t}^{0\mod 1}$ and,
after tedious but straightforward calculations, get
$$
\mu(R)\sim 0.
$$
Since $\mu(R)$ is quasimodular of weight four
and $\sim$ is equality modulo forms of weight
less than four, we get $\mu(R)=0$.
\end{proof}

\begin{remark}
There is a natural projection map from the space
of quasimodular forms of weight four to the space
of modular forms of weight four, which sends 
all forms divisible by $E_2$ to zero. So one can 
compose $\mu$ with this projection and have a 
map $\mu_1$ to the space of modular forms
of weight four.
\end{remark}

\begin{proposition}\label{minzero}
Map $\mu$ sends $M_4(l)_-$ to zero.
\end{proposition}

\begin{proof}
The statement immediately 
follows from the symmetry properties of $r$, $s$ and $t$.
\end{proof}

\begin{proposition}\label{immu}
The image of $S_4(l)_+$ under $\mu$ is 
a subspace of ${\mathcal S}_4(l)$
which is the linear span of $W(s_a,s_b)$ with
$gcd(a,b,l)=1$.
\end{proposition}

\begin{proof}
By Proposition \ref{xyspan}, $\mu(S_4(l))$
is spanned by $\mu(xy(a,b))=W(s_a,s_b)$ for all
$gcd(a,b,l)=1$. By Proposition \ref{minzero}, 
$\mu(S_4(l))_+=\mu(S_4(l))$.
\end{proof}

\section{The composition map}\label{sec.comp}
In this section we calculate the composition of 
the duality map $PD:S_4(l)_-^*\to S_4(l)_+$ of Section
\ref{sec.PD} and
the Wronskian map $\mu$ of Section \ref{sec.mu}. Our arguments
are purely elementary. The result of this calculation will
be used in the next section.

We need to introduce some additional notation.
For any $\phi\in S_4(l)_-^*$ we will set $\phi(P(x,y)(u,v)_-)=0$
if $gcd(u,v,l)>1$. We will also use the
notation $\sim$, namely, $f\sim g$ would mean
that $f-g$ is a linear combination of quasimodular forms
of weight at most two and the quasimodular forms
of weight three that are derivatives of $s_u(\tau)$.
Finally, for every $n>0$ we introduce 
the set $H(n)$ of fourtuples of integers $(a,b,c,d)$
that satisfy $ad-bc =n, a>b\geq 0,d>c\geq 0$.
\begin{proposition}\label{composition}
For any $\phi\in S_4(l)_-^*$ there holds
$$\mu\circ PD(\phi)\sim
\sum_{n>0}q^n\sum_{H(n)}\phi ((ax+by)(cx+dy)(c,d)_-).
$$ 
\end{proposition}

\begin{proof}
We will use the notations $A_{k_1,k_2}$ and $I(n)$
from Section \ref{sec.mu}. We use the last identity in
the proof of Proposition \ref{pm} to get 
$$
12\mu\circ PD(\phi)=
\mu(\sum_{u,v\in \ZZ/l\ZZ}(
\phi( (y-x)^2(-v,u+v)_-)x^2(u,v)_+
$$$$-2\phi( y(y-x)(-v,u+v)_-)xy(u,v)_+
+\phi( y^2(-v,u+v)_-) y^2(u,v)_+))
$$
$$
\sim
2\sum_{n>0}q^n
\phi\Big(\sum_{I(n)} 
\Big(
(-2m_1k_2(y-x)^2
-2(m_1k_1-m_2k_2)y(y-x)
$$$$+2m_2k_1y^2)
(-k_2,k_1+k_2)_-
+
(-2m_1k_2(y-x)^2
+2(m_1k_1-m_2k_2)y(y-x)
$$$$+2m_2k_1y^2
)
(-k_2,-k_1+k_2)_-)
\Big)
\Big)
$$$$
-\frac 1l\sum_{u,v\in \ZZ/l\ZZ}\phi((y-x)^2(-v,u+v)_-
-y^2(-u,v+u)_-) 
q\frac {\partial r_v}{\partial q} 
$$$$-\sum_{u\in\ZZ/l\ZZ}\phi((y-x)^2(0,u)_--y^2(-u,u)_-) 
q\frac {\partial t_u}{\partial q} .
$$

Let us simplify the last two lines of the above equation.
We use
$y^2(-u,v+u)_-= -x^2(-v,-u)_--(x-y)^2(v+u,-v)_-$,
$y^2(-u,u)_-= -x^2(0,-u)_--(x-y)^2(u,0)_-$,
$x^2(v,u)_-+x^2(v,-u)_-=y^2(v,u)_-+y^2(v,-u)_-=0$
and $x^2(0,u)_-=y^2(0,u)_-=0$ and $xy(0,u)_-=xy(u,0)_-$
to rewrite them as 
$$\frac 4l\sum_{u,v\in \ZZ/l\ZZ}\phi(xy(v,u)_-)
q\frac {\partial r_v}{\partial q}
+4\sum_u \phi(xy(u,0)_-)q\frac {\partial t_u}{\partial q}.
$$
\draft{
$$
-\frac 1l\sum_{u,v\in \ZZ/l\ZZ}\phi((y-x)^2(-v,u+v)_--y^2(-u,v+u)_-) 
q\frac {\partial r_v}{\partial q} 
$$$$=-\frac 1l\sum_{u,v\in \ZZ/l\ZZ}\phi((y-x)^2(-v,u+v)_-
+x^2(-v,-u)_-+(x-y)^2(v+u,-v)_-)q\frac {\partial r_v}{\partial q} 
$$$$=\frac 4l\sum_{u,v\in \ZZ/l\ZZ}\phi(xy(v,u)_-)
q\frac {\partial r_v}{\partial q}.$$
where we also use
$x^2(v,u)_-+x^2(v,-u)_-=y^2(v,u)_-+y^2(v,-u)_-=0$.
We have $y^2(-u,u)_-= -x^2(0,-u)_--(x-y)^2(u,0)_-$
so the last line is 
$$-\sum_{u}
\phi((y-x)^2(0,u)_-+x^2(0,-u)_-+(x-y)^2(u,0)_-) 
q\frac {\partial t_u}{\partial q}
$$
$$
=4\sum_u \phi(xy(u,0)_-)q\frac {\partial t_u}{\partial q}
$$
where we used
$x^2(0,u)_-=y^2(0,u)_-=0$ and $xy(0,u)_-=xy(u,0)_-$.
}

To handle the sum over $I(n)$, for each $n$ we observe
that $I(n)$ can be embedded into the disjoint union of two copies
of $H(n)$ in two different ways as follows.
The subset of $I(n)$ with $m_1\geq m_2$ 
can be identified with the subset of $H(n)$ with $c>0$
via $(m_1,k_1,m_2,k_2) = (a,d-c,a-b,c)$.
The subset of $I(n)$ with $m_1< m_2$ can
be identified with the subset of $H(n)$ with $bc>0$
via $(m_1,k_1,m_2,k_2) = (a-b,c,a,d-c)$. This 
describes the first 
embedding of  $I(n)$ into the disjoint union of two copies 
of $H(n)$. The second embedding is obtained by 
comparing $k_i$. Namely,
the subset of $I(n)$ with $k_1> k_2$ can
be identified with the subset of $H(n)$ with $bc>0$
via $(m_1,k_1,m_2,k_2)=(a-b,d,b,d-c)$, and
the subset of $I(n)$ with $k_1\leq k_2$ can be identified with
the subset of $H(n)$ with $b>0$ via
$(m_1,k_1,m_2,k_2)=(b,d-c,a-b,d)$.
We will use these embeddings in order to rewrite the 
above sum over $I(n)$ in terms of $H(n)$ as follows.
For the terms with $(-k_2,k_1+k_2)_-$ we will 
use the first embedding, and for the terms with 
$(-k_2,-k_1+k_2)_-$ we will use the second one.
After some straightforward simplifications, we get
$$
12\mu\circ PD(\phi)
\sim
4\sum_{n>0}q^n
\phi\Big(\sum_{H(n),bc>0}
\Big(
(-acx^2 +(ad+bc)xy -bdy^2)(-c,d)_-
$$$$
+
((-ad-bc+ac+bd)x^2+(ad+bc-2bd)xy+bdy^2)(c-d,d)_-
$$$$
+
((-ad-bc+ac+bd)x^2+(ad+bc-2ac)xy+acy^2)(c-d,-c)_-
$$$$
+
(-bdx^2+(ad+bc)xy-acy^2)(-d,c)_-
\Big)
$$$$
+\sum_{H(n),b=0,c>0} 
(-acx^2 +(ad+bc)xy -bdy^2)(-c,d)_-
$$
$$
+\sum_{H(n),b>0,c=0} 
((-bd)x^2+(ad+bc)xy+(-ac)y^2)(-d,c)_-
$$
$$
+\sum_{d\vert n}
\Big(\frac {2n^2}{d}xy(0,d)_-+\frac {2nd}{l}\sum_{u\in \ZZ/l\ZZ}
xy(d,u)_- 
\Big)\Big)
$$
$$
=
4\sum_{n>0}q^n
\phi\Big(\sum_{H(n),bc>0}
\Big(
(2acx^2 +2(ad+bc)xy+2bdy^2)(c,d)_-
$$$$+
((-ad-bc+ac+bd)x^2+(ad+bc-2bd)xy+bdy^2)
(c-d,d)_-
$$$$
+
((-ad-bc+ac+bd)x^2+(ad+bc-2ac)xy+acy^2)
)
(c-d,-c)_-
\Big)
$$
$$
+\sum_{H(n),bc=0} 
(acx^2 +(ad+bc)xy+bdy^2)(c,d)_-
-\sum_{d\vert n}nxy(0,d)_-
$$
$$
+\sum_{d\vert n}
\Big(\frac {2n^2}{d}xy(0,d)_-+\frac {2nd}{l}\sum_{u\in \ZZ/l\ZZ}
xy(d,u)_- 
\Big)\Big).
$$
We used various symmetries of $(u,v)_-$ to derive the 
last identity. A fortunate observation allows one to simplify
the second and third lines of the last formula. Indeed,
the relations on modular symbols imply
$$
((-ad-bc+ac+bd)x^2+(ad+bc-2bd)xy+bdy^2)(c-d,d)_- 
$$
$$
+(acx^2+(ad+bc-2ac)xy
+(-ad-bc+ac+bd)y^2)(-c,c-d)_- 
$$
$$
=(acx^2+(ad+bc)xy+bdy^2)(c,d)_-.
$$
Then one gets
$$
12\mu\circ PD(\phi)
\sim 
4\sum_{n>0}q^n
\phi\Big(\sum_{H(n)}
(3acx^2 +3(ad+bc)xy+3bdy^2)(c,d)_-
$$
$$
-2\sum_{H(n),bc=0} 
(acx^2 +(ad+bc)xy+bdy^2)(c,d)_-
-\sum_{d\vert n}nxy(0,d)_-
$$
$$
+\sum_{d\vert n}
\Big(\frac {2n^2}{d}xy(0,d)_-+\frac {2nd}{l}\sum_{u\in \ZZ/l\ZZ}
xy(d,u)_- 
\Big)\Big)
$$
$$
=4\sum_{n>0}q^n
\phi\Big(\sum_{H(n)}
(3acx^2 +3(ad+bc)xy+3bdy^2)(c,d)_-
$$
$$
-2\sum_{d\vert n,d>c>0} 
\frac nd(cx^2 +dxy)(c,d)_-
-\sum_{d\vert n}nxy(0,d)_-
$$$$
+\sum_{d\vert n}\frac {2nd}{l}
\sum_{u\in \ZZ/l\ZZ}xy(u,d)_- \Big).
$$

Using calculations similar to that of Section \ref{sec.mu},
we can rewrite the last two lines in terms of fractional parts as
\draft{
the sum over $d\vert n$
and $u\in\ZZ/l\ZZ$ of 
$$
n\Big(-\frac {2}d\sum_{0<c<d}(cx^2+dxy)\delta_c^{u\ml}
-xy\delta_0^{u\ml} + \frac {2d}lxy
\Big)
(u,d)_- 
$$
$$
=n\Big(-\frac {2}d\sum_{0<c\leq d}(cx^2+dxy)\delta_c^{u\ml}
+2 (x^2+xy)\delta_d^{u\ml}
-xy\delta_0^{u\ml} + \frac {2d}lxy
\Big)
(u,d)_- 
$$
...Intermediates
$$\sum_{0<c=rl+u\leq d,-u/l<r\leq (d-u)/l}
((rl+u)x^2+dxy) 
$$
$$= \sum_{r=[-\frac ul]+1}^{[\frac {d-u}l]}
((rl+u)x^2+dxy) 
$$
$$
= \frac 12 ([\frac {d-u}l]-[-\frac ul])
(((\frac {d-2u}l-\{\frac {d-u}l\}-\{-\frac ul\}+1)l+2u)x^2+2dxy)
$$
$$
= \frac 12 (\frac dl-\{\frac {d-u}l\}+\{-\frac ul\})
((d+l-l\{\frac {d-u}l\}-l\{-\frac ul\})x^2+2dxy)
$$
...End Intermediates
$$
=n\Big(-\frac {1}d(\frac dl-\{\frac {d-u}l\}+\{-\frac ul\})
(d+l-l\{\frac {d-u}l\}-l\{-\frac ul\})x^2
$$$$-2(\frac dl-\{\frac {d-u}l\}+\{-\frac ul\})xy
$$$$
+2 (x^2+xy)\delta_d^{u\ml}
-xy\delta_0^{u\ml} + \frac {2d}lxy
\Big)
(u,d)_- 
$$
$$
=n\Big(-\frac {1}d(\frac dl-\{\frac {d-u}l\}+\{-\frac ul\})
(d+l-l\{\frac {d-u}l\}-l\{-\frac ul\})x^2
$$$$+
(2\{\frac {d-u}l\}-2\{-\frac ul\}+2\delta_d^{u\ml}-\delta_0^{u\ml})xy
$$$$
+2 x^2\delta_d^{u\ml}
 \Big)
(u,d)_- 
$$
$$
=n\Big((-\frac {l}d(\frac dl-\{\frac {d-u}l\}+\{-\frac ul\})
(\frac dl+1-\{\frac {d-u}l\}-\{-\frac ul\})+2\delta_d^{u\ml})x^2
$$$$+
(-2\{\frac{u-d}l\}+1-\{-\frac {u}l\}+\{\frac ul\})xy
 \Big)
(u,d)_- 
$$
$$
=n\Big((
\frac {l}d\{\frac {d-u}l\}-\frac ld\{-\frac ul\}+\frac {l-d}{l}
-2\{\frac{u-d}l\}
-\frac {l}d\{\frac {d-u}l\}^2+\frac {l}d\{-\frac ul\}^2)
$$$$
+
(-2\{\frac{u-d}l\}+1-\{-\frac {u}l\}+\{\frac ul\})xy
 \Big)
(u,d)_- 
$$
}
$$
S=\phi\Big(
4\sum_{n>0}nq^n\sum_{d\vert n}\sum_{u\in\ZZ/l\ZZ}\Big((
\frac {l}d(\{\frac {d-u}l\}-\{\frac {d-u}l\}^2)
-\frac ld(\{-\frac ul\}-\{-\frac ul\}^2)
$$$$+\frac {l-d}{l}
-2\{\frac{u-d}l\})x^2
+
(-2\{\frac{u-d}l\}+1-\{-\frac {u}l\}+\{\frac ul\})xy
 \Big)
(u,d)_-\Big).
$$
We observe that for any $t$ there holds 
$\{t\}-\{t\}^2 = \{-t\}-\{-t\}^2$, and then use symmetries
of $P(x,y)(\pm u,d)_-$ to see that 
$$
S\sim\phi\Big(
4\sum_{n>0}nq^n\sum_{d\vert n}\sum_{u\in\ZZ/l\ZZ}\Big(
-2\{\frac{u-d}l\}x^2+
(-2\{\frac{u-d}l\}+1)xy
 \Big)
(u,d)_-\Big)
$$
$$
=\phi
\Big(4
\sum_{n>0}nq^n\sum_{d\vert n}\sum_{u\in\ZZ/l\ZZ}\Big(
(\{\frac{-u-d}l\}-\{\frac{u-d}l\})x^2+
$$$$(-\{\frac{u-d}l\} -\{\frac{-u-d}l\}+1)xy
 \Big)
(u,d)_-\Big)
$$
$$
\sim\phi\Big(
4\sum_{n>0}nq^n\sum_{d\vert n}\sum_{u\in\ZZ/l\ZZ}\Big(
(\{1-\delta_{d}^{-u\ml})x^2+\delta_{d}^{u\ml} xy
 \Big)
(u,d)_-\Big)
$$
$$
=\phi\Big(
4\sum_{n>0}nq^n\sum_{d\vert n}(x^2+xy)(d,d)_-)\Big)\sim 0.
$$
This finishes the proof.
\end{proof}

\section{Relation to Hecke eigenforms of rank zero}
In this section we prove our main result that relates
Wronskians of weight one Eisenstein series and
the Hecke eigenforms of weight four with nonzero
central value of $L$-function.

Let $T_n$ denote the Hecke operators for $\Gamma_1(l)$
and let $L(f,s)$ denote the Hecke $L$-function. We will normalize
it so the central value is $L(f,2)=\int_{0}^{\ii\infty}f(\tau)\tau \,d\tau$.
We say that a weight four Hecke eigenform $f$ has analytic
rank zero if $L(f,2)\neq 0$.
\begin{definition}\label{defrho}
Let $f\in{\mathcal S}_4(l)$ 
be a weight four  cusp form for $\Gamma_1(l)$.
Define $\rho(f)=\sum_{n>0} L(T_nf,2)q^n$.
\end{definition}

The following statements are analogous to
the weight two calculation of \cite{vanish}.
\begin{proposition}\label{likevanish}
Definition \ref{defrho}
gives a linear map $\rho:{\mathcal S}_4(l)
\to {\mathcal S}_4(l)$,
which commutes with $\Gamma_0(l)/\Gamma_1(l)$-action.
The image of $\rho$ contains all newforms $f$
with $L(f,2)\neq 0$, and is contained in the span of
all Atkin-Lehner lifts of all Hecke eigenforms $f$ of analytic
rank zero.
\end{proposition}

\begin{proof}
The arguments of \cite[Propositions 4.3 and 4.5]{vanish}
apply to weight four case without any serious changes.
\end{proof}

Similar to  \cite{vanish}, the key idea of this paper
is to relate the map $\rho$ to the map $\mu$ of 
Section \ref{sec.mu}. 
\begin{proposition}
The map $\rho$ is the composition of the maps
$$
{\mathcal S}_4(l)
\stackrel {Int} \to (S_4(l)_-)^*
\stackrel {PD}\to 
S_4(l)_+\stackrel{\mu}\to {\mathcal S}_4(l)
$$
where $Int$ is induced by the integration
pairing of ${\mathcal S}_4(l)$ and $S_4(l)_-$,
the $PD$ is the Poincar\'e duality map of Section \ref{sec.PD},
and $\mu$ is the Wronskian map of Section \ref{sec.mu}.
\end{proposition}

\begin{proof}
We denote by $\langle,\rangle$ the integration pairing
between ${\mathcal S}_4(l)$ and $S_4(l)_-$.
For a given $f\in {\mathcal S}_4(l)$ we calculate
$$
\rho(f) = \sum_{n>0} L(T_nf,2)q^n=\sum_{n>0} 
\langle T_nf,xy(0,1)_-\rangle q^n.
$$
By \cite[Theorem 2 and Proposition 10]{Merel}, 
$$\langle T_nf,xy(0,1)_-\rangle=
\langle f, T_nxy(0,1)_-\rangle =
\langle f,\sum_{H(n)} (ax+by)(cx+dy)(c,d)_-\rangle
$$
which leads to 
$$
\rho(f)=\sum_{n>0}q^n
\langle f,\sum_{H(n)} (ax+by)(cx+dy)(c,d)_-\rangle
$$
Proposition \ref{composition} now shows
$\mu\circ PD\circ Int(f)\sim \rho(f)$ and since both sides
are quasimodular forms of weight four, the claim follows.
\end{proof}

\begin{corollary}\label{imrho}
The image of $\rho$ equals the linear span 
of $W(s_a,s_b)$ for $gcd(a,b,l)=1$.
\end{corollary}

\begin{proof}
Recall that $Int$ and $PD$ are isomorphisms, by 
\cite{Merel} and Corollary \ref{surj} respectively.
Then Proposition \ref{immu} finishes the proof.
\end{proof}

We are now ready to formulate our main result.
\begin{theorem}\label{main}
For arbitrary $l>1$ the span of Hecke eigenforms
of weight four and analytic rank zero is equal to the
span of the Wronskians $W(s_a,s_b)$ for all
$a,b\in \ZZ/l\ZZ$.
\end{theorem}

\begin{proof}
In one direction, consider $f=W(s_a,s_b)$. If
$gcd(a,b,l)=d$, then Corollary \ref{imrho} applied
to $\frac ld$ shows
that $f$ is in $\rho({\mathcal S}_4(\frac ld))$.
Indeed, $f$ is, up to a nonzero multiple, 
the $d$-lift of 
$W(s_{\frac ad,\frac ld},s_{\frac bd,\frac ld})$
where the second subscript in $s$ is used to
indicate the level. By Proposition \ref{likevanish},
$W(s_{\frac ad,\frac ld},s_{\frac bd,\frac ld})$
lies in the linear span of eigenforms of analytic rank zero,
hence $f$ does as well.

To prove the opposite inclusion, 
it is enough to show that for any $d\vert l$
and any newform $g(\tau)\in {\mathcal S}_4(\frac ld)$
of analytic rank zero, its lift 
$g(k\tau)\in {\mathcal S}_4(l)$ lies in the span
of $W(s_a,s_b)$ for any $k\vert d$.
By Proposition \ref{likevanish}, $g\in\rho({\mathcal S}_4(\frac ld))$.
Then by Corollary \ref{imrho}, 
$g$ is a linear combination of Wronskians of Eisenstein
series $s_{i,\frac ld}$ of level $\frac ld$. 
Then $g(k\tau)$ is a linear combination 
of Wronskians of $s$-series of level $\frac {kl}d$, since 
$s_{i,\frac ld}(k\tau)=s_{ki,\frac {kl}d}(\tau)$.
Finally, $s$-series of level $\frac {kl}d$ are sums
of $s_a$ of level $l$, which shows that 
$g(k\tau)$ lies in the span of the Wronskians
$W(s_a,s_b)$, as claimed.
\end{proof}

\begin{corollary}
The span of Hecke eigenforms of weight four and analytic
rank zero for the group $\Gamma_0(l)$ coincides with the
span of
$$
\sum_{j\in(\ZZ/l\ZZ)^*} W(s_{aj},s_{bj})
$$
for all $a,b\in\ZZ/l\ZZ$.
\end{corollary}

\begin{proof}
Use the formulas for the action of $\Gamma_0(l)$
on $s_a$ from  \cite{highweight}.
\end{proof}

\end{document}